\documentclass[11pt]{article}
\usepackage{anysize}\marginsize{3.5cm}{3.5cm}{2cm}{3cm}
\usepackage{amssymb}
\newtheorem{satz}{Satz}[section]

\newtheorem{corollary}[satz]{Corollary}

\newtheorem{example}[satz]{Example}

\newtheorem{proposition}[satz]{Proposition}

\newtheorem{remarks}[satz]{Remarks}
\newtheorem{theorem}[satz]{Theorem}

\newtheorem{thevarthm}[satz]{\varthmname}

\newenvironment{varthm*}[1]{\trivlist\item[]{\bf #1.}\it}{\endtrivlist}

\renewcommand\ge{\geqslant}  
\renewcommand\le{\leqslant}  
\renewcommand\geq{\geqslant}  
\renewcommand\leq{\leqslant}  
\renewcommand\epsilon{\varepsilon}
\renewcommand\phi{\varphi}

\renewcommand\O{{\cal O}}

\newcommand\be{\begin{eqnarray*}}
\newcommand\ee{\end{eqnarray*}}

\newcommand\eqnref[1]{(\ref{#1})}
\newcommand\eqdef{=_{\rm def}}
\newcommand\eps{\varepsilon}
\renewcommand\P{\mathbb P}
\newcommand\N{\mathbb N}

\newcommand\inparen[1]{\textnormal{(}{#1}\textnormal{)}}
\newcommand\mapto{\longrightarrow}
\newcommand\mult{{\rm mult}}

\begin{document}

\title{Seshadri constants on surfaces of general type}
\author{Thomas Bauer and Tomasz Szemberg}
\date{December 13, 2007}
\maketitle
\thispagestyle{empty}

\begin{abstract}
   We study Seshadri constants of the canonical bundle on minimal
   surfaces of general type. First, we prove that if the
   Seshadri constant $\eps(K_X,x)$ is between 0 and 1, then it is
   of the form $(m-1)/m$ for some integer $m\ge 2$.
   Secondly, we study values of $\eps(K_X,x)$ for a very general point
   $x$ and show that small values of the Seshadri constant are accounted for
   by the geometry of $X$.
\end{abstract}

\section*{Introduction}

\begin{plaintheoremnumbers}

   Given a smooth projective variety $X$ and a nef line bundle $L$ on
   $X$, Demailly defines the Seshadri constant of $L$
   at a point $x\in X$ as the real number
   $$
      \eps(L,x)\eqdef\inf_C\frac{L\cdot C}{\mult_xC} \ ,
   $$
   where the infimum is taken over all irreducible
   curves passing through $x$
   (see \cite{Dem92} and \cite[Chapt.~5]{PAG}).
   If $L$ is ample, then
   $\eps(L,x)>0$ for all points $x\in X$.
   When $X$ is a surface, then by a result of Ein and Lazarsfeld 
   \cite{EL93} one has for ample $L$
   $$
      \eps(L,x)\ge 1 \quad\mbox{for all except perhaps countably many points
      $x\in X$.}
      \eqno(*)
   $$

   In the present paper we study the Seshadri constants of the
   canonical bundle~$K_X$ on minimal surfaces of general type
   (i.e., we study $\eps(K_X,x)$ on 
   surfaces $X$ of Kodaira dimension~2 that do
   not contain any $(-1)$-curves, or, equivalently, on surfaces
   $X$ whose
   $K_X$ is big and nef).
   Motivated by the result $(*)$ of Ein-Lazarsfeld
   we first ask for
   the potential values below 1 that $\eps(K_X,x)$ might have. We
   show:

\begin{theorem}\label{arbitrary}
   Let $X$ be a smooth projective surface such that the canonical
   divisor $K_X$ is big and nef.
   Let $x$ be any point on $X$.
\begin{itemize}
\item[\rm(a)]
   One has $\eps(K_X,x)=0$ if and only if $x$ lies on one of the
   finitely many $(-2)$-curves on $X$.
\item[\rm(b)]
   If $0<\eps(K_X,x)<1$, then there is an integer $m\ge 2$ such
   that
   $$
      \eps(K_X,x)=\frac{m-1}m \ ,
   $$
   and there is an irreducible curve $C\subset X$ such that $\mult_x(C)=m$ and
   $K_X\cdot C=m-1$. \inparen{In other words, the curve $C$ computes
   the Seshadri constant of $K_X$ at $x$.}
\item[\rm(c)]
   If $0<\eps(K_X,x)<1$ and $K_X^2\ge 2$, then either
   \begin{itemize}
   \item[\rm(i)]
      $\eps(K_X,x)=\frac12$ and $x$ is the double point of
      an irreducible
      curve $C$ with arithmetic genus $p_a(C)=1$ and $K_X\cdot
      C=1$, or
   \item[\rm(ii)]
      $\eps(K_X,x)=\frac 23$ and $x$ is a triple point of 
      an irreducible 
      curve $C$ with arithmetic genus $p_a(C)=3$ and $K_X\cdot C=2$.
   \end{itemize}
\item[\rm(d)]
   If $0<\eps(K_X,x)<1$ and $K_X^2\ge 3$, then only case {\rm(c)(i)} is possible.
\end{itemize}
\end{theorem}

   It is well known that the bicanonical system $|2K_X|$ is base
    point free on almost all surfaces of general type. For such
   surfaces one gets easily the lower bound $\eps(K_X,x)\geq 1/2$
   for all $x$ away of the contracted locus.
   In general however one knows only that $|4K_X|$ is base point
   free, which gives a lower bound of $1/4$. Theorem~\ref{arbitrary} shows
   in particular that
   one has $\eps(K_X,x)\ge 1/2$ in all cases. Moreover
   this bound turns out to be sharp, i.e.,
   there are examples of surfaces $X$ and points $x$
   such that $\eps(K_X,x)=1/2$ (see Example~\ref{example half}).
   We do not know whether all values $(m-1)/m$ for arbitrary $m\ge
   2$ actually occur.
   As part (c) of Theorem \ref{arbitrary} shows, however, 
   values $(m-1)/m$ with $m\ge 4$ can occur
   only in the case $K_X^2=1$.
   We will show in Example~\ref{doubleconic} 
   that curves as in (c)(i) actually exist on surfaces with
   arbitrarily large degree of the canonical bundle.
   In other words, one cannot strengthen the result by imposing
   higher bounds on $K_X^2$.
   It would be interesting to know
   whether curves as in (c)(ii)
   exist.

   We consider next Seshadri constants of the canonical
   bundle at a very general point. It is known
   that the function $x\mapsto\eps(L,x)$ is
   lower semi-continuous in the topology on $X$, whose closed
   sets are the countable unions of subvarieties
   (see \cite{Ogu02}). In particular, it assumes its
   maximal value for $x$
   very general, i.e., away of an at most countable union of proper Zariski
   closed subsets of $X$. This maximal value will be denoted by $\eps(L,1)$.
   We show:

\begin{theorem}\label{general}
   Let $X$ be a smooth projective surface such that $K_X$ is big
   and nef. If $K_X^2\ge 2$, then
   $$
      \eps(K_X,1)>1 \ .
   $$
\end{theorem}

   Note that it may well happen that $\eps(K_X,x)=1$ for
   infinitely many $x\in X$ even when $K_X^2\ge 2$:
   Consider for instance
   a smooth quintic surface $X\subset\P^3$ containing a
   line $\ell$. Then $\eps(K_X,x)=1$ for all $x\in\ell$. However,
   one has $\eps(K_X,1)>1$ by \cite[Theorem~2.1(a)]{Bau}.
   Theorem~\ref{general} states that the same conclusion holds in
   general
   as soon as $K_X^2\ge 2$.

   If $K_X^2$ is even larger, then further
   geometric statements are possible. We show:

\begin{theorem}\label{thm fibration}
   Let $X$ be a smooth projective surface such that $K_X$ is big
   and nef. If $K_X^2\ge 6$, then
   \begin{itemize}
   \item[\rm(a)]
      $\eps(K_X,1)\ge 2$,
   \item[\rm(b)]
      $\eps(K_X,1)=2$ if and only if $X$ admits a genus 2
      fibration $X\to B$ over a smooth curve $B$.
   \end{itemize}
\end{theorem}

   A somewhat more general statement is given in Propositions
   \ref{zwei} and \ref{fibration}.
\end{plaintheoremnumbers}

\paragraph*{Acknowledgement.}
   This paper was written during the authors' visits in Krak\'ow and
   Marburg.
   The first author was partially supported by DFG grant BA 1559/4-3,
   the second author by
   KBN grant 1 P03 A 008 28. We are grateful to the referee for
   pointing out omissions in the first version of Theorem 1(c).
   Further, we would like to thank
   Fabrizio Catanese for sharing with us the idea of using conic
   bundles in order to construct example \ref{doubleconic}.

\section{Seshadri constants at arbitrary points}

   In this section we prove Theorem~\ref{arbitrary}.
   We will need the following fact:

\begin{proposition}\label{computing-curve}
   Let $L$ be a nef and big
   line bundle on a smooth projective surface, and let $x$ be
   a point
   such that $\eps(L,x)<\sqrt{L^2}$. Then there is an irreducible curve $C$
   such that
   $$
      \eps(L,x)=\frac{L\cdot C}{\mult_x(C)} \ .
   $$
\end{proposition}

   The point here is that $\eps(L,x)$ is in fact
   computed by a curve
   rather than being approximated by a sequence of curves.

\begin{proof}
   We fix a real number $\xi$ such that
   $1<\xi<\frac{\sqrt{L^2}}{\eps(L,x)}$.
   By definition of $\eps(L,x)$ there is in any event a sequence
   $(C_n)_{n\in\N}$ of irreducible
   curves such that $L\cdot C_n/\mult_x(C)$ converges
   to $\eps(L,x)$ from above. In particular $L\cdot C_n/\mult_x(C_n)<\frac{1}{\xi}\cdot\sqrt{L^2}$
   for $n\gg 0$.

   On the other hand the asymptotic Riemann-Roch theorem implies that for $k\gg 0$
   there are divisors $D\in|kL|$ such that the quotient
   $L\cdot D /\mult_x D$ is arbitrarily close to $\sqrt{L^2}$,
   in particular less than $\xi\cdot\sqrt{L^2}$.
   Fixing such a value $k$,
   \cite[Lemma 5.2]{Bau}
   implies that
   for all sufficiently large $n$, the curve $C_n$ is a component
   of $D$. Therefore in the sequence $(C_n)$ there are only
   finitely many distinct curves, and this implies the assertion.
\end{proof}

\begin{proofof}{Theorem~\ref{arbitrary}}
   (a)
   Suppose first that $C$ is a rational $(-2)$-curve on $X$. Then
   $K_X\cdot C=0$ by the adjunction formula, so that clearly
   $$
      \eps(K_X,x)=0
   $$
   for any point $x\in C$. Conversely, suppose that
   $\eps(K_X,x)=0$ for some point $x\in X$.
   Thanks to Proposition~\ref{computing-curve} 
   we have
   $K_X\cdot C=0$
   for some irreducible curve $C$ on $X$.
   As $K_X^2>0$, we get $C^2<0$ from the index theorem. The
   adjunction formula then tells us that $C$ is a $(-2)$-curve.

   (b)
   Since by assumption $0<\eps(K_X,x)<1$, there is in any event
   an integer
   $m\ge 2$ such that
   \begin{equation}\label{interval}
      \frac{m-2}{m-1}<\eps(K_X,x)\le\frac{m-1}m \ .
   \end{equation}
   We will show that then necessarily $\eps(K_X,x)=
   (m-1)/m$.
   To this end we start by making
   use of a result from \cite{Bau} (cf.\ Proof of \cite[Theorem
   3.1]{Bau}):
   \begin{itemize}\item[]
      If $L$ is a nef and big line bundle,
      and if $\sigma$ is a real number such that
      $\sigma L-K_X$ is nef, then for
      any irreducible curve $C$ the
      multiplicity $\ell=\mult_x(C)$ at a given point $x\in X$
      is bounded as
      \begin{equation}\label{m_bound}
         \ell\le\frac12+\sqrt{\frac{(L\cdot C)^2}{L^2}+\sigma L\cdot C+\frac94} \ .
      \end{equation}
   \end{itemize}
   (Note that in \cite{Bau} the line bundle $L$ is assumed to be
   ample. The argument proving
   \eqnref{m_bound} still holds, however,
   when $L$ is merely nef and big.)
   In our situation we take $L=K_X$ and $\sigma=1$, and we have
   by Proposition~\ref{computing-curve} an irreducible curve
   $C$ such that
   $$
      \eps(K_X,x)=\frac{K_X\cdot C}{\ell} \ .
   $$
   Note that since $\eps(K_X,x)$ is less than $1$ we have $\ell\geq
   2$. Moreover if $\ell=2$, then the only possible value of the
   Seshadri constant is $1/2$. From now on we assume therefore that
   $\ell\geq 3$.

   Writing $d=K_X\cdot C$, the inequality \eqnref{m_bound}
   implies
   then
   \begin{equation}\label{eps_bound}
      \frac{K_X\cdot C}{\ell}\ge\frac{d}{
         \frac12+\sqrt{\frac{d^2}{K_X^2}+d+\frac94}} \ .
   \end{equation}
   The expression on the right hand side of \eqnref{eps_bound} is
   increasing both as a function of $d$ and as a function
   of $K_X^2$.
   Now, a calculation shows that
   if $d\ge m$, then the right hand side of
   \eqnref{eps_bound} is bigger than $(m-1)/m$.
   So our assumption \eqnref{interval} implies that
   $$
      d\le m-1 \ .
   $$
   We will now show that in fact
   $d= m-1$.
   Indeed, the Seshadri constant $\eps(K_X,x)$ is the fraction
   $$
      \eps(K_X,x)=\frac d\ell \ ,
   $$
   but among all fractions that are  smaller than 1 and that
   have numerator $d$, the fraction $d/(d+1)$ is the biggest. So
   $$
      \eps(K_X,x)\le\frac d{d+1} \ .
   $$
   But if $d\le m-2$, then this implies
   $$
      \eps(K_X,x)\le\frac{m-2}{m-1} \ ,
   $$
   giving a contradiction with \eqnref{interval}.
   We have thus established the equality
   $$
      d=m-1 \ .
   $$
   From the index theorem we obtain
   \begin{equation}\label{hit}
      K_X^2\,C^2\le (K_X\cdot C)^2=d^2 \ .
   \end{equation}
   Note that this implies that in any event $C^2\le d^2$.
   The adjunction formula then gives
   $$
      p_a(C)=1+\frac 12\, C\cdot(C+K_X)\le 1+\frac{d(d+1)}2 \ .
   $$
   As $C$ has a point of
   multiplicity $\ell$,
   the arithmetic genus of its normalization
   drops by at least $({\ell\atop 2})$. Since $\ell\geq 3$, we get therefore
   $$
      \ell\le d+1 = m \ .
   $$
   So we conclude that
   $$
      \eps(K_X,x)=\frac d\ell\ge\frac dm=\frac{m-1}m \ ,
   $$
   and this completes the proof of (b).

   (c)
   The argument is parallel to that of part~(b) up to the
   inequality \eqnref{hit}. We keep the notations from there.
   Suppose that $0<\eps(K_X,x)<1$ and $K_X^2\ge 2$.
   We see from \eqnref{hit} that
   $C^2\le \frac12 d^2$. 
   If we have $\ell\ge 4$, then we obtain
   \begin{equation}\label{ell-inequality}
      \ell<d+1=m \ ,
   \end{equation}
   and therefore
   $$
      \eps(K_X,x)=\frac d\ell>\frac dm=\frac{m-1}m \ ,
   $$
   in contradiction with \eqnref{interval}.

   So it remains to study the cases $\ell=2$ and $\ell=3$. 
   If $\ell=2$, 
   then necessarily
   $d=K_X\cdot C=1$, so that $\eps(K_X,x)=1/2$.
   The assumption $K_X^2\geq 2$ together with \eqnref{hit} imply that
   $C^2<0$. The adjunction formula and the fact that $C$ has a
   double point imply then that $p_a(C)=1$.
   
   If $\ell=3$, then we arrive at the situation $d=K_X\cdot C=2$,
   $C^2=2$, and hence $p_a(C)=3$.

   (d)
   Keeping notations from the preceding case, the assumption
   $K_X^2\ge 3$ implies now via \eqnref{hit} the inequality
   $C^2\le\frac13 d^2$, so that we arrive for $\ell\ge 3$ again
   at \eqnref{ell-inequality}, which gives a contradiction as
   before.
\end{proofof}

   We now describe an example of a surface $X$ of general type,
   where one has
   $$
      K_X^2=1 \qquad\mbox{and}\qquad \eps(K_X,x)=\frac 12 \qquad\mbox{for some $x\in X$.}
   $$
\begin{example}\label{example half}\rm
   Consider a
   general surface $X$ of degree 10
   in weighted projective space $\P(1,1,2,5)$.
   Then $X$ is smooth, $K_X=\O_X(10-1-1-2-5)=\O_X(1)$
   by adjunction. In particular, $K_X$ is ample, $K_X^2=1$, and $h^0(K_X)=2$,
   which corresponds to the first two variables of weight $1$
   (see \cite[p. 311]{Ste81} for details).
   The canonical pencil consists of curves of arithmetic genus $2$
   and has exactly one base point in which all canonical curves
   are smooth. Blowing up this point $\sigma:Y\mapto X$ we get a
   genus $2$ fibration $f:Y\mapto \P^1$ over the exceptional curve.
   The mapping $f$ is a relatively minimal semistable family of
   curves. In fact there are neither multiple nor reducible fibers
   possible, as
   $K_X^2=1$ and $K_X$ is ample. So the singular fibers, if any,
   are irreducible curves with double points.
   If all fibers were smooth genus $2$ curves, then the
   topological Euler characteristic of $Y$ would be
   $c_2(Y)=2\cdot (-2)=-4$.
   Since $\chi(Y)=\chi(X)=3$, this would contradict the Noether formula.
   This shows that there exists a singular
   canonical curve $D\in|K_X|$. In the singular point of $D$ one
   has the Seshadri quotient $\frac12$.
\end{example}

   We now give an example showing that
   curves as in part (c)(i) and (d) of Theorem~\ref{arbitrary} occur on surfaces with
   arbitrarily large degree of the canonical bundle.

\begin{example}\label{doubleconic}\rm
\newcommand\calo{{\cal O}}
\newcommand\calc{{\cal C}}
\newcommand\cald{{\cal D}}
   In the product $\P^2\times\P^2$ we consider
   a nontrivial conic bundle $\calc_0\mapto \P^2$ over
   the first factor with a smooth discriminant curve $\Delta\subset\P^2$.
   This can be obtained explicitly as a general divisor of bidegree $(1,2)$.
   Let $B\subset\P^2$ be a smooth plane curve of degree $d\geq 4$
   intersecting $\Delta$ transversally. We restrict $\calc_0$ to a conic
   bundle $f:\calc\mapto B$ over $B$. We fix a point $b_0\in
   B\cap\Delta$. The fiber of $f$ over $b_0$ consists of two lines
   $L_0$ and $L_1$ meeting in a point $P$. It is easy to find a smooth
   cubic $D$ in the plane spanned by $L_0$ and $L_1$
   not passing through $P$ and such that it
   intersects $L_0$ in three distinct points and such that it cuts out on $L_1$ 
   a divisor of the form $2Q+R$. We extend $D$ to a hypersurface in
   $\P^2\times\P^2$ simply taking the product
   $\cald=\P^2\times D$. If $D$ is sufficiently
   general, then the intersection curve $\Gamma=\calc\cap\cald$ is
   smooth. Taking the double covering $\sigma:X\mapto \calc$
   branched over $\Gamma$ we obtain a smooth surface $X$ with a
   genus $2$ fibration $\alpha:X\mapto B$. By the subadditivity
   of Kodaira dimensions $X$ is of general type. This construction 
   may be considered as reversing the
   procedure described before Theorem 4.13 in
   \cite{CatPig06}. The fiber over $b_0$
   consists of two reduced and irreducible curves $C_0$ and $C_1$
   lying over $L_0$ and $L_1$ respectively. The curve $C_0$ is a smooth
   elliptic curve, while $C_1$ has a double point over $Q$ and
   arithmetic genus $1$. Since $(C_0+C_1)^2=0$ and $(C_0+C_1)\cdot K_X=2$,
   we see by adjunction that $C_1^2=-1$ and $K_X\cdot C_1=1$.
   Hence $C_1$ is a curve as in part (c) of Theorem
   \ref{arbitrary}. The formula in \cite[Theorem 4.13]{CatPig06}
   shows that $K_X^2$ grows linearly with the genus
   of $B$, so that
   taking the curve $B$ of high enough degree we can make $K_X^2$ arbitrarily large.
\end{example}

\section{Seshadri constants at very general points}

   In this section we prove Theorem \ref{general}, which is
   stated in the
   introduction.
   We start by recalling the following
   result obtained by Syzdek and the second author \cite{SyzSze}.
   The result is stated in \cite{SyzSze} for ample line bundles,
   but the proof works verbatim for big and nef ones.

\begin{theorem}\label{ss}
   Let $L$ be a big and nef line bundle on a smooth projective
   surface $X$. Assume that
   $$\eps(L,1)<\sqrt{\frac79L^2}.$$
   Then $X$ is either a smooth cubic in $\P^3$
   or $X$ is fibred by curves computing $\eps(L,x)$.
\end{theorem}

\begin{proofof}{Theorem~\ref{general}}
   We claim first that there are at most finitely many reduced and irreducible
   curves $C\subset X$ with $K_X\cdot C\le 1$.
   To see this,
   consider first the case that $K_X\cdot C=0$. Then $C$ is a
   $(-2)$-curve, and we know that there are only finitely many of
   them on $X$. Next, suppose $K_X\cdot C=1$. Then the index
   theorem gives $K_X^2 C^2\le (K_X\cdot C)^2=1$, which implies
   $C^2\le 0$. From the genus formula
   $p_a(C)=1+\frac12(C^2+K_X\cdot C)$ we see that $C^2$ is an odd
   number. So we have $C^2\le -1$, and
   therefore $C$ is the only irreducible curve in its
   numerical equivalence class. The claim now follows
   from the fact that
   there are only finitely many classes in the Neron-Sev\'eri
   group of $X$ that have degree $1$ with respect to $K_X$.

   Now assume to the contrary that $\eps(K_X,1)=1$.
   Then the numerical assumptions of Theorem~\ref{ss}
   are satisfied for the line bundle
   $L=K_X$. Since a cubic in $\P^3$ is not of general type,
   there must be a fibration of curves computing $\eps(K_X,x)$.
   Since any curve $C$ in the fibration is smooth in its general point
   it must be $K_X\cdot C=1$ which contradicts above reasoning
   on the number of such curves.
\end{proofof}

   The proof shows that in the situation of the theorem
   one has in fact the lower bound
   $\eps(K_X,x)\ge\frac 13\sqrt{14}$. It is
   unlikely, however, that this
   particular bound is sharp.

   Theorem~\ref{general} gives
   in particular
   the following interesting characterization of the
   situations in which equality holds in statement~$(*)$ at the beginning
   of the introduction.

\begin{corollary}
   Let $X$ a smooth projective surface such that $K_X$ is big and
   nef. Then $K_X^2=1$ if and only if $\eps(K_X,1)=1$.
\end{corollary}

\begin{remarks}\rm
   (i) The corollary should be seen in light of the fact that in
   general it is very
   well possible to have ample line bundles $L$
   on smooth projective surfaces such that $\eps(L,x)=1$ for
   very general $x$, while $L^2$ can be arbitrarily large.
   Consider for instance a product $X=C\times D$ of two smooth irreducible curves,
   and denote by a slight abuse of notation the fibers of both
   projections again by $D$ and $C$.
   The line bundles $L_m=mC+D$ are ample and we have $L_m\cdot C=1$, so that
   in any event $\eps(L_m,x)\le 1$ for every point $x\in X$.
   One has in fact $\eps(L_m,x)=1$, which can be seen as follows:
   If $F$ is any irreducible
   curve different from the fibers of the projections
   with $x\in F$, then we may take a 
   fiber $D'$ of the first projection with $x\in D'$, and we have
   $$
      L_m\cdot F \ge D'\cdot F \ge \mult_x(D')\cdot\mult_x(F)
      \ge\mult_x(F)
   $$
   which implies $\eps(L_m,x)\ge 1$.
   So $\eps(L_m,x)=1$, but 
   on the other hand $L_m^2=2m$
   is unbounded.

   (ii) Consider for a moment a minimal surface $X$
   of general type
   such that $p_g=0$.  One knows then that the bicanonical
   system $|2K_X|$ is composed with a pencil if and only if
   $K_X^2=1$ (see \cite[Theorem~3.1]{MenLop04}).
   The corollary shows that this geometric
   condition is also encoded in
   Seshadri constants through the condition
   $\eps(K_X,1)=1$.
\end{remarks}

   In the spirit of Theorem~\ref{general} we now show the following result.

\begin{proposition}\label{zwei}
   Let $X$ be a smooth projective surface such that $K_X$ is big and nef.
   If $K_X^2\geq 6$, then either
   \begin{itemize}
   \item[\rm(a)] $\eps(K_X,1)>2$, or
   \item[\rm(b)] $\eps(K_X,1)=2$, and there exists a pencil of curves of genus $2$
   computing $\eps(K_X,x)$ for $x$ very general.
   \end{itemize}
\end{proposition}
\begin{proof}
   If $\eps(K_X,1)> 2$, then we are done. So let us assume that
   $\eps(K_X,1)\leq 2$. Then the numerical assumptions of Theorem~\ref{ss}
   are satisfied, and hence there is a fibration of curves computing $\eps(K_X,1)$.
   Since a curve in the fibration is smooth in its general point $x$,
   $\eps(K_X,x)$ is an integer. By Theorem~\ref{general} we have
   then
   $\eps(K_X,1)=2$, and $\eps(K_X,x)$ is
   computed by a curve $C$ satisfying $K_X\cdot C=2$ and $C^2=0$.
   Hence $C$ is a member of a pencil of curves of genus $2$.
\end{proof}

   The presence of a genus $2$ fibration on a surface of general type
   is the typical obstacle in
   the positivity of the canonical bundle and the associated maps. Hence
   it is not very surprising that it implies an upper bound on
   $\eps(K_X,1)$, which can be regarded as a converse of
   Proposition
   \ref{zwei}(b).

\begin{proposition}\label{fibration}
   Let $X$ be a smooth minimal surface of general type such that
   there is a genus $2$ fibration $f:X\to B$ over a smooth curve
   $B$. Then
   $$\eps(K_X,1)\leq 2 \ ,$$
   and if $K_X^2\geq 4$, then actually
   $$\eps(K_X,1)=2\ .$$
\end{proposition}

\begin{proof}
   If $K_X^2\leq 4$, then the assertion is clear, because in any event
   one has $\eps(K_X,x)\le\sqrt{K_X^2}$ for all $x\in X$ by
   Kleiman's theorem (cf.~\cite[Proposition~5.1.9]{PAG}).
   So we may assume that $K_X^2\geq 5$. Let $F$ be a generic fiber of
   $f$. As $K_X$ restricts to the canonical bundle on $F$, we
   have
   $K_X\cdot F=2$.
   In order to conclude we need to show that for a general point $x$
   the quotient $\frac{K_X\cdot C}{\mult_xC}$ is greater or equal
   $2$ for all irreducible curves passing through the point $x$.
   Since on a surface of general type there is no elliptic or rational curve
   passing through a general point, we have by adjunction $K_X\cdot
   C\geq 2$. Suppose
   for a contradiction
   that $\frac{K_X\cdot C}{\mult_xC}<2$ holds.
   We must then in particular have $\mult_xC\geq 2$.
   As $\frac{K_X\cdot C}{\mult_xC}<\sqrt{K_X^2}$,
   the argument from the proof of the main theorem in \cite{EL93}
   shows that there exists
   a nontrivial family of pointed curves $(C,x)$ with
   multiplicity $m\geq 2$ at $x$ such that
   $\frac{K_X\cdot C}{m}<2$. Then
   $K_X\cdot C\leq 2m-1$ and, by \cite[Corollary 1.2]{EL93},
   $C^2\geq m(m-1)$ . Combining this
   with the index theorem we have
   $$5m(m-1)\leq K_X^2 C^2\leq (K_X\cdot C)^2\leq (2m-1)^2 \ ,$$
   which is impossible for $m\geq 2$.
\end{proof}

   Propositions \ref{zwei} und \ref{fibration} imply
   Theorem~\ref{thm fibration} from the introduction.

\bigskip
\small
   Tho\-mas Bau\-er,
   Fach\-be\-reich Ma\-the\-ma\-tik und In\-for\-ma\-tik,
   Philipps-Uni\-ver\-si\-t\"at Mar\-burg,
   Hans-Meer\-wein-Stra{\ss}e,
   D-35032~Mar\-burg, Germany.

\nopagebreak
   E-mail: \texttt{tbauer@mathematik.uni-marburg.de}

\bigskip
   Tomasz Szemberg,
   Instytut Matematyki AP,
   PL-30-084 Krak\'ow, Poland

\nopagebreak
   E-mail: \texttt{szemberg@ap.krakow.pl}
\end{document}